# Probability matching priors for some parameters of the bivariate normal distribution

Malay Ghosh[*,1], Upasana Santra[1] and Dalho Kim[2]

*University of Florida and Kyungpook National University*

**Abstract:** This paper develops some objective priors for certain parameters of the bivariate normal distribution. The parameters considered are the regression coefficient, the generalized variance, and the ratio of the conditional variance of one variable given the other to the marginal variance of the other variable. The criterion used is the asymptotic matching of coverage probabilities of Bayesian credible intervals with the corresponding frequentist coverage probabilities. The paper uses various matching criteria, namely, quantile matching, matching of distribution functions, highest posterior density matching, and matching via inversion of test statistics. One particular prior is found which meets all the matching criteria individually for all the parameters of interest.

## Contents

## 1. Introduction

Bayesian methods are becoming increasingly popular in the theory and practice of statistics. It is needless to say that other than the likelihood, the key component in any Bayesian analysis is the selection of priors. With sufficient information from past experience, expert opinion or previously collected data, subjective priors are ideal, and indeed should be used for inferential purposes. However, often even without adequate prior information, one can use Bayesian techniques efficiently with some

[*]Supported in part by the NSF Grant SES-06-31426.
[1]Department of Statistics, University of Florida, Gainesville, Florida 32611-8545, USA, e-mail: ghoshm@stat.ufl.edu
[2]Department of Statistics, Kyungpook National University, Taegu, Korea

"default" or "objective" priors. Not too surprisingly, the catalog of such priors, over the years, has also become prohibitively large, and one needs to set some specific "objectivity" criterion for the solution of such priors.

One such criterion which has found some appeal to both frequentists and Bayesians is the so-called "probability matching" criterion. Simply put, this amounts to the requirement that the coverage probability of a Bayesian credible region is asymptotically equivalent to the coverage probability of the frequentist confidence region up to a certain order. An excellent monograph on this topic is due to Datta and Mukerjee [5] which provides a thorough and comprehensive discussion of various probability matching criteria. Other review papers include [6], [7] and [9].

Again, as one might expect, there are several probability matching criteria. The matching is accomplished through either (a) posterior quantiles, (b) distribution functions, (c) highest posterior density (HPD) regions, or (d) inversion of certain test statistics. However, priors based on (a), (b), (c), or (d) need not always be identical. Specifically, it may so happen that there does not exist any prior satisfying all four criteria.

In this article, we consider the bivariate normal distribution where the parameters of interest are the (i) regression coefficient, (ii) the generalized variance, i.e. the determinant of the variance-covariance matrix, and (iii) the ratio of the conditional variance of one variable given the other divided by the marginal variance of the other variable. We have been able to find a prior which meets all four matching criteria for every one of these parameters.

As is well known, matching priors are obtained by solving certain differential equations. In Section 2 of this paper, we have introduced an orthogonal reparameterization ([3], [8]) of the bivariate normal parameters which greatly simplifies these equations, resulting thereby in easily accessible solutions. Here, we have also introduced a quantile matching prior which works for all the parameters given in (i)–(iii). Section 3 establishes the distribution matching property of the prior found in Section 2, once again for all three parameters of interest. Section 4 establishes the HPD matching property of such priors, while Section 5 confirms matching by inversion of the likelihood ratio statistic. The propriety of the posteriors and some numerical computations are given in Section 6, while some final remarks are made in Section 7.

## 2. The orthogonal reparameterization

Let $(X_{1i}, X_{2i}), i = 1, \ldots, n$, be independent and identically distributed random variables having a bivariate normal distribution with means $\mu_1$ and $\mu_2$, variances $\sigma_1^2(>0)$ and $\sigma_2^2(>0)$, and correlation coefficient $\rho(|\rho|<1)$. We use the transformation (see, e.g., [1], [4] or [10])

$$(2.1) \qquad \beta = \rho\sigma_2/\sigma_1,\ \theta = \sigma_1\sigma_2(1-\rho^2)^{1/2} \text{ and } \eta = \sigma_2(1-\rho^2)^{1/2}/\sigma_1.$$

With this reparameterization, the bivariate normal distribution can be rewritten as
(2.2)
$$f(X_1, X_2) = (2\pi\theta)^{-1} \exp\left\{-\frac{1}{2}\left\{\frac{(X_2 - \mu_2 - \beta(X_1 - \mu_1))^2}{\theta\eta} + \frac{\eta(X_1-\mu_1)^2}{\theta}\right\}\right\}$$

where $\beta$ is the regression coefficient of $X_2$ on $X_1$, $\theta^2 = \sigma_1^2\sigma_2^2(1-\rho^2)$ is the determinant of the variance-covariance matrix, and $\eta^2 = V(X_2|X_1)/V(X_1)$.

With the above reparameterization, the Fisher Information matrix reduces to

$$\mathbf{I}(\mu_1, \mu_2, \beta, \theta, \eta) = \begin{pmatrix} \mathbf{A} & \mathbf{0} \\ \mathbf{0} & \text{Diag}(\eta^{-2}, \theta^{-2}, \eta^{-2}) \end{pmatrix}, \quad (2.3)$$

where

$$\mathbf{A} = \begin{pmatrix} \frac{\beta^2}{\theta\eta} + \frac{\eta}{\theta} & -\frac{\beta}{\theta\eta} \\ -\frac{\beta}{\theta\eta} & \frac{1}{\theta\eta} \end{pmatrix}.$$

This establishes immediately the mutual orthogonality of $(\mu_1, \mu_2)$, $\beta$, $\theta$ and $\eta$ in the sense of [3] and [8]. Such orthogonality is often referred to as "Fisher Orthogonality."

Since the parameters of interest are orthogonal to $(\mu_1, \mu_2)$, and it is customary to use a uniform $(\Re^2)$ prior on $(\mu_1, \mu_2)$, we shall consider only priors of the form

$$\pi_0(\mu_1, \mu_2, \beta, \theta, \eta) \propto \pi(\beta, \theta, \eta), \quad (2.4)$$

and find $\pi$ such that the matching criteria given in (a)-(d) are all satisfied for $\beta$, $\theta$ and $\eta$, each individually. This we are going to explore in the next three sections.

Here we also state a lemma which is used repeatedly in the sequel. The proof is based on the independence of $X_2 - \mu_2 - \beta(X_1 - \mu_1)$ and $X_1 - \mu_1$ along with the fact that $X_2 - \mu_2 - \beta(X_1 - \mu_1) \sim N(0, \theta\eta)$ and $X_1 - \mu_1 \sim N(0, \theta/\eta)$.

**Lemma 2.1.** *For the bivariate normal density given in (2.2),*

$$E(\partial \log f/\partial \beta)^3 = 0, \ E[(\partial \log f/\partial \beta)(\partial^2 \log f/\partial \beta^2)] = 0; \quad (2.5)$$

$$E(\partial^3 \log f/\partial \beta^3) = 0, \ E(\partial^3 \log f/\partial \beta^2 \partial \theta) = (\theta\eta^2)^{-1}, \ E(\partial^3 \log f/\partial \beta^2 \partial \eta) = \eta^{-3}; \quad (2.6)$$

$$E(\partial^3 \log f/\partial \beta \partial \theta^2) = 0, E(\partial^3 \log f/\partial \beta \partial \eta^2) = 0; \quad (2.7)$$

$$E(\partial \log f/\partial \theta)^3 = 2/\theta^3, \ E[(\partial \log f/\partial \theta)(\partial^2 \log f/\partial \theta^2)] = -2/\theta^3; \quad (2.8)$$

$$E(\partial^3 \log f/\partial \theta^3) = 4/\theta^3, \ E(\partial^3 \log f/\partial \theta^2 \partial \eta) = 0, \ E(\partial^3 \log f/\partial \theta \partial \eta^2) = (\theta\eta^2)^{-1}; \quad (2.9)$$

$$\begin{aligned} E(\partial \log f/\partial \eta)^3 &= 0, \\ E[(\partial \log f/\partial \eta)(\partial^2 \log f/\partial \eta^2)] &= -\eta^{-3}, E(\partial^3 \log f/\partial \eta^3) = 3\eta^{-3}. \end{aligned} \quad (2.10)$$

We consider in this section quantile matching priors, i.e., priors which ensure approximate frequentist validity of one-sided Bayesian credible intervals based on posterior quantiles of a one-dimensional interest parameter. The pioneering research in this area is due to Welch and Peers [14] and Peers [11], while the recent stimulus comes from Stein [12] and Tibshirani [13]. Specifically, one considers here priors $\pi(\cdot)$ for which the relation

$$P_\theta\{\theta_1 \leq \theta_1^{1-\alpha}(\pi, X)\} = 1 - \alpha + o(n^{-\frac{r}{2}}) \quad (2.11)$$

holds for $r = 1$ or 2 and for each $\alpha$ $(0 < \alpha < 1)$. In the above, $\theta = (\theta_1, \ldots, \theta_p)^T$ is the unknown parameter, $\theta_1$ is the one-dimensional parameter of interest, $\theta_1^{1-\alpha}(\pi, X)$ is the $(1 - \alpha)$th posterior quantile of $\theta_1$ based on the prior $\pi$ and data $X =$ $(X_1, \ldots, X_n)^T$, while $P(\cdot|\theta)$ denotes the conditional probability given $\theta$, the usual frequentist probability. A prior satisfying (2.11) with $r = 1$ is called a *first order* probability matching prior, while one with $r = 2$ is called a *second order* probability matching prior. Clearly, second order probability matching priors constitute a subclass of first order probability matching priors.

A second order quantile matching prior which works for $\beta$, $\theta$ and $\eta$ is given by $\pi(\mu_1, \mu_2, \beta, \theta, \eta) \propto (\theta\eta)^{-1}$. Back to the original parameterization, $(\mu_1, \mu_2, \sigma_1, \sigma_2, \rho)$, this reduces to the prior

$$(2.12) \qquad \pi(\mu_1, \mu_2, \sigma_1, \sigma_2, \rho) = \sigma_1^{-2}(1 - \rho^2)^{-1}.$$

This prior has been identified in [2] as the right-Haar prior as well as the one-at-a-time reference prior for any arbitrary ordering of these parameters. Indeed, this prior provides exact rather than just asymptotic matching for a variety of parameters of interest including the ones considered here. Moreover, as shown in this paper, when $\beta$ is the parameter of interest, any prior of the form $\sigma_1^{-(3-a)}(1-\rho^2)^{-1}$, $a > 0$, for $(\mu_1, \mu_2, \sigma_1, \sigma_2, \rho)$ achieves exact matching, while when $\theta$ is the parameter of interest, both the priors $\sigma_1^{-2}(1-\rho^2)^{-1}$ and $\sigma_1^{-1}\sigma_2^{-1}(1-\rho^2)^{-3/2}$ achieve exact matching.

However, Berger and Sun [1] have not explored other matching criteria. While the quantile matching property is quite desirable for construction of one-sided credible intervals, the HPD matching or matching via inversion of test statistics seems more appropriate for the construction of two-sided credible intervals. It is also true that priors which satisfy a particular matching criterion, for example, the quantile matching criterion, satisfies other matching criteria. As mentioned in the introduction, we will explore various other matching criteria in the subsequent sections. In this way, we will lend further justification to one of the matching priors of [2] – the one that we have developed here, as well.

## 3. Matching via distribution functions

In this section, we target priors $\pi$ which achieve matching via distribution functions of some standardized variables. More specifically, when $\theta_1$ is the parameter of interest, $(\theta_2, \ldots, \theta_p)^T$ is the vector of nuisance parameters, $\hat{\theta}_1$ is the MLE of $\theta_1$ with $n^{-\frac{1}{2}}I^{11}$ $(I = ((I_{jj'})), I^{-1} = ((I^{jj'})))$ as its asymptotic variance, we consider the random variable $y = \sqrt{n}(\theta_1 - \hat{\theta}_1)/(I^{11})^{1/2}$. Specifically, if $P^\pi$ denotes the posterior of $y$ given the data X, what we want to achieve is the asymptotic matching

$$(3.1) \qquad E[P^\pi(y \leq w|X)|\theta] = P(y \leq w|\theta) + o(n^{-1}).$$

Under orthogonality of $\theta_1$ with $(\theta_2, \ldots, \theta_p)$, it follows from (3.2.5) to (3.2.7) of [5] that such a prior $\pi$ is of the form $I_{11}^{1/2}g(\theta_2, \ldots, \theta_p)$, where in addition one needs to satisfy the two differential equations

(3.2)
$$A_1 = \frac{\partial^2}{\partial\theta_1^2}\left(I^{11}\pi(\theta)\right) - 2\frac{\partial}{\partial\theta_1}\left(I^{11}\pi(\theta)\right) - \sum_{s=2}^{p}\sum_{v=2}^{p}\frac{\partial}{\partial\theta_s}\left\{E\left(\frac{\partial^3 \log f}{\partial\theta_1^2 \partial\theta_s}\right)I^{11}I^{sv}\pi(\theta)\right\}$$
$$- \sum_{s=2}^{p}\sum_{v=2}^{p}\frac{\partial}{\partial\theta_1}\left\{E\left(\frac{\partial^3 \log f}{\partial\theta_1 \partial\theta_s \partial\theta_v}\right)I^{11}I^{sv}\pi(\theta)\right\} = 0$$

and

$$(3.3) \quad A_2 = \sum_{s=2}^{p}\sum_{v=2}^{p} \frac{\partial}{\partial \theta_s}\left\{E\left(\frac{\partial^3 \log f}{\partial \theta_1^2 \partial \theta_s}\right) I^{11} I^{sv} \pi(\theta)\right\} = 0.$$

In the given problem, when $\beta$ is the parameter of interest, (3.3) reduces to

$$(3.4) \quad \frac{\partial}{\partial \beta}\left(\eta^2\left\{\theta^2 E\left(\frac{\partial^3 \log f}{\partial \beta^2 \partial \theta}\right) + \eta^2 E\left(\frac{\partial^3 \log f}{\partial \beta^2 \partial \eta}\right)\right\}\pi(\beta,\theta,\eta)\right) = 0.$$

By (2.6) of Lemma 2.1, (3.4) simplifies to $\partial/\partial\beta\{(\theta+\eta)\pi(\beta,\theta,\eta)\} = 0$ which holds trivially for any prior $\pi(\beta,\theta,\eta)$ which does not depend on $\beta$, including the prior $\pi(\beta,\theta,\eta) \propto (\theta\eta)^{-1}$, the one found in Section 2. Again, with $\beta$ as the parameter of interest, for any prior $\pi(\beta,\theta,\eta)$ which does not depend on $\beta$, (3.2) simplifies to $\partial/\partial\theta\{(\theta\eta^2)^{-1}\eta^2\theta^2\pi(\beta,\theta,\eta)\} + \partial/\partial\eta\{\eta^{-3}\eta^2\eta^2\pi(\beta,\theta,\eta)\} = 0$, i.e.,

$$\frac{\partial}{\partial \theta}[\theta\,\pi(\beta,\theta,\eta)] + \frac{\partial}{\partial \eta}[\eta\,\pi(\beta,\theta,\eta)] = 0.$$

Once again $\pi(\beta,\theta,\eta) \propto (\theta\eta)^{-1}$ will satisfy (3.3). To verify that a prior $\pi(\beta,\theta,\eta)$ which does not depend on $\beta$ satisfies (3.2) with $\theta_1 = \beta$, we need to verify that

$$(3.5) \quad \frac{\partial}{\partial \beta}\left(\eta^2\left\{E\left(\frac{\partial^3 \log f}{\partial \beta \partial \theta^2}\right)\theta^2 + E\left(\frac{\partial^3 \log f}{\partial \beta \partial \eta^2}\right)\eta^2\right\}\right) = 0,$$

which reduces to

$$\eta^2\theta^2 \frac{\partial}{\partial \beta}E\left(\frac{\partial^3 \log f}{\partial \beta \partial \theta^2}\right) + \eta^4 \frac{\partial}{\partial \beta}E\left(\frac{\partial^3 \log f}{\partial \beta \partial \eta^2}\right) = 0.$$

From (2.7) of Lemma 2.1, $\mathrm{E}(\partial^3 \log f/\partial\beta\partial\theta^2) = \mathrm{E}(\partial^3 \log f/\partial\beta\partial\eta^2) = 0$, so that (3.2) holds trivially for such a prior. Hence matching via distributions is achieved with any prior of the form $\pi(\mu_1,\mu_2,\beta,\theta,\eta) \propto h(\theta,\eta)$, and in particular $h(\theta,\eta) \propto (\theta\eta)^{-1}$

Next, when $\theta$ is the parameter of interest, to find a prior satisfying (3.1) one needs to first solve

$$(3.6) \quad \begin{aligned}&\frac{\partial^2}{\partial \theta^2}\left(\theta^2\pi(\cdot)\right) - 2\frac{\partial}{\partial \theta}\left(\theta^2 \frac{\partial \pi(\cdot)}{\partial \theta}\right) - \frac{\partial}{\partial \theta}\left(\theta^4 \mathrm{E}\left(\frac{\partial^3 \log f}{\partial \theta^3}\right)\pi(\cdot)\right) \\ &\quad - \frac{\partial}{\partial \theta}\left(\theta^4 \mathrm{E}\left(\frac{\partial^3 \log f}{\partial \theta^3}\right)\pi(\cdot)\right) = 0.\end{aligned}$$

Hence, (3.6) simplifies to

$$\frac{\partial^2}{\partial \theta^2}\left(\theta^2\pi(\cdot)\right) - 2\frac{\partial}{\partial \theta}\left(\theta^2 \frac{\partial \pi(\cdot)}{\partial \theta}\right) - 12\frac{\partial}{\partial \theta}\left(\theta\pi(\cdot)\right) = 0.$$

Any prior $\pi(\cdot) \propto \theta^{-1}g(\beta,\eta)$ will satisfy this equation. Such a prior also satisfies

$$(3.7) \quad \frac{\partial}{\partial \theta}\left(\theta^4 E\left(\frac{\partial^3 \log f}{\partial \theta^3}\right)\pi(\cdot)\right) = 0.$$

Finally when $\eta$ is the parameter of interest, again, for finding a prior satisfying (3.1), one needs to solve

$$
\begin{aligned}
(3.8) \quad & \frac{\partial^2}{\partial \eta^2}\left\{\eta^2 \pi(\cdot)\right\} - 2\frac{\partial}{\partial \eta}\left\{\eta^2 \frac{\partial \pi(\cdot)}{\partial \eta}\right\} - \frac{\partial}{\partial \theta}\left\{E\Big(\frac{\partial^3 \log f}{\partial \theta \partial \eta^2}\Big)\eta^2 \theta^2 \pi\right\} \\
& - \frac{\partial}{\partial \eta}\left\{E\Big(\frac{\partial^3 \log f}{\partial \eta \partial \beta^2}\Big)\eta^2 \eta^2 \pi\right\} = 0
\end{aligned}
$$

and

$$
(3.9) \quad \frac{\partial}{\partial \eta}\Big(\eta^4 E\Big(\frac{\partial^3 \log f}{\partial \eta^3}\Big)\pi(\cdot)\Big) = 0.
$$

Again, by Lemma 2.1, the prior $\pi(\beta, \theta, \eta) = \theta^{-1}\eta^{-1}$ satisfies the matching property.

## 4. Highest posterior density (HPD) matching priors

We now turn our attention to HPD matching priors for each one of the parameters $\beta, \theta$ and $\eta$. In general, if $\tilde{\theta}$ is the parameter (real or vector-valued) of interest, then a HPD region is of the form $\{\tilde{\theta} : \pi(\tilde{\theta}|X) \geq k\}$, where $\pi(\tilde{\theta}|X)$ is the posterior of $\tilde{\theta}$ under the prior $\pi$ and data X. We will consider priors which ensure that HPD regions with credibility level $1 - \alpha$ also have asymptotically the same frequentist coverage probability, the error of approximation being $o(n^{-1})$.

We first consider a HPD region for $\beta$. In view of the orthogonality result derived in Section 2, such a prior $\pi_0(\mu_1, \mu_2, \beta, \theta, \eta) \propto \pi(\beta, \theta, \eta)$ must satisfy (see [5], (4.4.1))

$$
\begin{aligned}
(4.1) \quad & \frac{\partial}{\partial \theta}\left(\eta^2 \theta^2 \mathrm{E}\Big(\frac{\partial^3 \log f}{\partial \beta^2 \partial \theta}\Big)\pi\right) + \frac{\partial}{\partial \eta}\left(\eta^4 \mathrm{E}\Big(\frac{\partial^3 \log f}{\partial \beta^2 \partial \eta}\Big)\pi\right) \\
& + \frac{\partial}{\partial \beta}\left(\eta^4 \mathrm{E}\Big(\frac{\partial^3 \log f}{\partial \beta^3}\Big)\pi\right) - \frac{\partial^2}{\partial \beta^2}\Big(\eta^2 \pi\Big) = 0.
\end{aligned}
$$

Again, by Lemma 2.1, (4.1) reduces to

$$
(4.2) \quad \frac{\partial}{\partial \theta}(\theta \pi(\cdot)) + \frac{\partial}{\partial \eta}(\eta \pi(\cdot)) - \frac{\partial^2}{\partial \beta^2}(\eta^2 \pi(\cdot)) = 0.
$$

Clearly the prior $\pi(\beta, \theta, \eta) \propto (\theta\eta)^{-1}$ satisfies (4.2).

Next consider $\theta$ as the parameter of interest. Now one needs to solve

$$
\begin{aligned}
(4.3) \quad & \frac{\partial}{\partial \beta}\left(\eta^2 \theta^2 \mathrm{E}\Big(\frac{\partial^3 \log f}{\partial \theta^2 \partial \beta}\Big)\pi\right) + \frac{\partial}{\partial \eta}\left(\theta^2 \eta^2 \mathrm{E}\Big(\frac{\partial^3 \log f}{\partial \theta^2 \partial \eta}\Big)\pi\right) \\
& + \frac{\partial}{\partial \theta}\left(\theta^4 \mathrm{E}\Big(\frac{\partial^3 \log f}{\partial \theta^3}\Big)\pi\right) - \frac{\partial^2}{\partial \theta^2}\Big(\theta^2 \pi\Big) = 0.
\end{aligned}
$$

By (2.8) and (2.9) of Lemma 2.1, (4.3) simplifies to

$$
(4.4) \quad -2\frac{\partial}{\partial \theta}(\theta \pi) - \frac{\partial^2}{\partial \theta^2}(\theta \pi) = 0,
$$

which is satisfied by the prior $\pi(\beta, \theta, \eta) \propto (\theta\eta)^{-1}$.

Finally, when $\eta$ is the parameter of interest, we need to solve

(4.5)
$$\frac{\partial}{\partial \beta}\left(\eta^4 E(\frac{\partial^3 \log f}{\partial \eta^2 \partial \beta})\pi\right) + \frac{\partial}{\partial \theta}\left(\eta^2 \theta^2 E(\frac{\partial^3 \log f}{\partial \eta^2 \partial \theta})\pi\right)$$
$$+ \frac{\partial}{\partial \eta}\left(\eta^4 E(\frac{\partial^3 log f}{\partial \eta^3})\pi\right) - \frac{\partial^2}{\partial \eta^2}\left(\eta^2 \pi\right) = 0.$$

From Lemma 2.1, (4.5) reduces to

(4.6)
$$\frac{\partial}{\partial \theta}(\theta \pi) + \frac{\partial}{\partial \eta}(\eta \pi) - \frac{\partial^2}{\partial \eta^2}(\eta^2 \pi) = 0.$$

Again $\pi(\beta, \theta, \eta) \propto (\theta \eta)^{-1}$ will do.

## 5. Matching priors via inversion of test statistics

One traditional way to derive frequentist confidence intervals is inversion of certain test statistics. The most popular such test is the likelihood ratio test. But tests based on Rao's score statistic or the Wald statistic are also of importance, and are first order equivalent (i.e., up to $o(n^{-1/2})$) to the likelihood ratio tests. We consider here only the likelihood ratio test.

We begin with the general case when $\theta_1$ is the parameter of interest, while $\theta_2, \ldots, \theta_p$ are the nuisance parameters. Let $\theta = (\theta_1, \ldots, \theta_p)$, and let $l(\theta)$ denote the usual log-likelihood. The corresponding profile log-likelihood for $\theta_1$ is given by $l^*(\theta_1) = l(\theta_1, \hat{\theta}_2(\theta_1), \ldots, \hat{\theta}_p(\theta_1))$, where $\hat{\theta}_j(\theta_1)$ is the MLE of $\theta_j$ given $\theta_1 (j = 2, \ldots, p)$. Then the likelihood ratio statistic for $\theta_1$ is given by

(5.1)
$$M_{LR}^*(\theta_1, X) = 2[l(\hat{\theta}) - l^*(\theta_1)].$$

Then from Yin and Ghosh [15] (also from [5], (5.2.18)), a likelihood ratio matching prior $\pi$ is obtained by solving

(5.2)
$$\sum_{s=2}^{p}\sum_{u=2}^{p}\frac{\partial}{\partial u}\left\{I_{11}^{-1}I^{su}E(\frac{\partial^3 \log f}{\partial \theta_1^2 \partial \theta_s})\pi\right\} + \frac{\partial}{\partial \theta_1}\left(I_{11}^{-1}\left\{\frac{\partial \pi}{\partial \theta_1} - \pi(\theta)(I_{11}^{-1}E(\frac{\partial \log f}{\partial \theta_1}\frac{\partial^2 \log f}{\partial \theta_1^2})\right.\right.$$
$$\left.\left. -\sum_{s=2}^{p}\sum_{u=2}^{p}I^{su}E(\frac{\partial^3 \log f}{\partial \theta_1 \partial \theta_u \partial \theta_s}))\right\}\right) = 0.$$

In the present case when $\beta$ is the parameter of interest, (5.2) reduces to

(5.3)
$$\frac{\partial}{\partial \theta}\left(\eta^2 \theta^2 (\theta \eta^2)^{-1}\pi\right) + \frac{\partial}{\partial \eta}\left(\eta^4 \eta^{-3}\pi\right) + \frac{\partial}{\partial \beta}\left(\eta^2 \left\{\frac{\partial \pi}{\partial \beta}\right.\right.$$
$$\left.\left. - \pi\left\{\eta^2 E((\frac{\partial \log f}{\partial \beta})(\frac{\partial^2 \log f}{\partial \beta^2})) - \theta^2 E(\frac{\partial^3 \log f}{\partial \beta \partial \theta^2}) - \eta^2 E(\frac{\partial^3 \log f}{\partial \beta \partial \eta^2})\right\}\right\}\right) = 0.$$

From Lemma 2.1, (5.3) reduces to

$$\frac{\partial}{\partial \theta}(\theta \pi) + \frac{\partial}{\partial \eta}(\eta \pi) + \eta^2 \frac{\partial}{\partial \beta}(\frac{\partial \pi}{\partial \beta}) = 0,$$

i.e.,
$$\frac{\partial}{\partial \theta}(\theta \pi) + \frac{\partial}{\partial \eta}(\eta \pi) + \eta^2 \frac{\partial^2 \pi}{\partial \beta^2} = 0.$$

Again $\pi \propto (\theta \eta)^{-1}$ provides a solution.

Next, if $\theta$ is the parameter of interest, the LR matching prior $\pi$ for $\theta$ is obtained by solving the differential equation

(5.4)
$$\frac{\partial}{\partial \beta}\left\{\eta^2 \theta^2 \cdot 0 \cdot \pi\right\} + \frac{\partial}{\partial \eta}\left\{\eta^2 \theta^2 \cdot 0 \cdot \pi\right\} + \frac{\partial}{\partial \theta}\left(\theta^2 \left\{\frac{\partial \pi}{\partial \theta}\right.\right.$$
$$\left.\left. - \pi \left\{\theta^2 E\left(\left(\frac{\partial \log f}{\partial \theta}\right)\left(\frac{\partial^2 \log f}{\partial \theta^2}\right)\right) - \eta^2\left(E\left(\frac{\partial^3 \log f}{\partial \beta^2 \partial \theta}\right) + E\left(\frac{\partial^3 \log f}{\partial \eta^2 \partial \theta}\right)\right)\right\}\right\}\right) = 0.$$

Again from Lemma 2.1, (5.4) reduces to
$$\frac{\partial}{\partial \theta}[\theta^2 \{\frac{\partial \pi}{\partial \theta} + \frac{2}{\theta}\pi + \pi \frac{2}{\theta}\}] = 0,$$

i.e.,
$$\frac{\partial}{\partial \theta}[\theta^2 \frac{\partial \pi}{\partial \theta} + 4\theta \pi] = 0$$

which holds for $\pi \propto (\theta \eta)^{-1}$.

Finally, when $\eta$ is the parameter of interest, the LR matching prior is obtained by solving

(5.5)
$$\frac{\partial}{\partial \beta}\left\{\eta^4 \eta^{-3} \cdot 0 \cdot \pi\right\} + \frac{\partial}{\partial \theta}\left\{\eta^2 \theta^2 \frac{1}{\theta \eta^2} \pi\right\}$$
$$+ \frac{\partial}{\partial \eta}\left(\eta^2 \left\{\frac{\partial \pi}{\partial \eta} - \pi \left\{\eta^2 E\left(\left(\frac{\partial \log f}{\partial \eta}\right)\left(\frac{\partial^2 \log f}{\partial \eta^2}\right)\right) - \eta^2 \frac{1}{\eta^3} - 0\right\}\right\}\right) = 0.$$

Once again, using Lemma 2.1, (5.5) reduces to
$$\frac{\partial}{\partial \theta}(\theta \pi) + \frac{\partial}{\partial \eta}[\eta^2(\frac{\partial \pi}{\partial \eta}) - \pi \eta^2(-\frac{2}{\eta})] = 0,$$

and the prior $\pi \propto (\theta \eta)^{-1}$ provides a solution.

## 6. Posteriors and numerically computed coverage

The prior $\pi(\mu_1, \mu_2, \beta, \theta, \eta) \propto (\theta \eta)^{-1}$ is improper. In this section, we write down the marginal posteriors for $\beta, \theta$ and $\eta$, and discuss methods for finding the HPD intervals for each one of these parameters. The analytical findings of Section 4 are strengthened with some numerical coverage probability computations.

The marginal posterior of $\beta$ is given by
$$\pi(\beta|\mathbf{X}_1, \mathbf{X}_2) \propto \int_o^\infty \eta^{n-2}\left(\eta^2 + \frac{S_{22} + \beta^2 S_{11} - 2\beta S_{12}}{S_{11}}\right)^{-(n-1)} d\eta.$$

Putting $\eta = z[S_{22} + \beta^2 S_{11} - 2\beta S_{12}/S_{11}]^{-1/2}$ in the above integral, one gets after simplification,

$$\pi(\beta|\mathbf{X}_1, \mathbf{X}_2) \propto \left(1 + \frac{(\beta - S_{12}/S_{11})^2}{S_{22.1}}\right)^{-\frac{n-1}{2}}, \tag{6.1}$$

where $S_{22.1} = S_{22} - S_{12}^2/S_{11}$. This posterior is a t-distribution with location parameter $S_{12}/S_{11}$, scale parameter $\{S_{22.1}/(n-2)\}^{1/2}$ and degrees of freedom $n-2$.

The posterior of $\theta$ is given by

$$\pi(\theta|\mathbf{X}_1, \mathbf{X}_2) \propto \theta^{-(n-1)} \exp(-S_{11}^{1/2} S_{22.1}^{1/2}/\theta) I_{[\theta > 0]}, \tag{6.2}$$

so that $\theta^{-1}$ has a Gamma distribution with shape parameter $n-2$ and scale parameter $(S_{11} S_{22.1})^{-\frac{1}{2}}$.

Finally, the marginal posterior of $\eta$ is given by

$$\begin{aligned}\pi(\eta|\mathbf{X}_1, \mathbf{X}_2) &\propto \eta^{-1/2} \left(\frac{S_{22.1}}{\eta} + \eta S_{11}\right)^{-(n-3/2)} \\ &\propto \eta^{n-2} \left(\eta^2 + \frac{S_{22.1}}{S_{11}}\right)^{-(n-3/2)}.\end{aligned} \tag{6.3}$$

The construction of HPD credible intervals is fairly simple. The posterior of $\beta$ being a univariate-$t$ (thus symmetric and unimodal), from (6.1), the $100(1-\alpha)\%$ HPD credible interval for $\beta$ is given by $S_{12}/S_{11} \pm \{S_{22.1}/(n-2)\}^{1/2} t_{n-2;\alpha/2}$, where $t_{n-2;\alpha/2}$ denotes the upper $100(\alpha/2)\%$ point of a Student's t-distribution with $n-2$ degrees of freedom.

Observing that the posterior of $\theta$ is log-concave, the $100(1-\alpha)\%$ region for $\theta$ is given by $[\theta_1, \theta_2]$, where $\theta_1$ and $\theta_2$ satisfy

$$\theta_1^{-(n-1)} \exp(-S_{11}^{1/2} S_{22.1}^{1/2}/\theta_1) = \theta_2^{-(n-1)} \exp(-S_{11}^{1/2} S_{22.1}^{1/2}/\theta_2) \tag{6.4}$$

and

$$\int_{\theta_1}^{\theta_2} \theta^{-(n-1)} \exp(-S_{11}^{1/2} S_{22.1}^{1/2}/\theta)(S_{11} S_{22.1})^{\frac{n-2}{2}} d\theta = 1 - \alpha. \tag{6.5}$$

It is important to note that if $w = \theta^{-1}$, then the posterior pdf of $w$ is given by

$$\pi(w|\mathbf{X}_1, \mathbf{X}_2) \propto w^{n-3} \exp(-w S_{11}^{1/2} S_{22.1}^{1/2}).$$

Noting the log-concavity of this pdf as well, the HPD region $[w_1, w_2]$ for $w$ is obtained by solving

$$w_1^{n-3} \exp(-w_1 S_{11}^{1/2} S_{22.1}^{1/2}) = w_2^{n-3} \exp(-w_2 S_{11}^{1/2} S_{22.1}^{1/2}) \tag{6.6}$$

and

$$\int_{w_1}^{w_2} \frac{w^{n-3}}{\Gamma(n-2)} \exp(-w S_{11}^{1/2} S_{22.1}^{1/2})(S_{11} S_{22.1})^{\frac{n-2}{2}} dw = 1 - \alpha. \tag{6.7}$$

Clearly the solution $[w_1, w_2]$ of (6.6) and (6.7) is different from the solution $[\theta_2^{-1}, \theta_1^{-1}]$ of (6.4) and (6.5).

TABLE 1
*Frequentist coverage probabilities of 95% HPD intervals
for $\beta$, $\theta$ and $\eta$ when $\sigma_1^2 = 1$ and $\sigma_2^2 = 1$*

| $\rho$ | $n$ | $\beta$ | $\theta$ | $\eta$ |
|---|---|---|---|---|
| 0.25 | 4 | 0.952 | 0.947 | 0.949 |
| | 8 | 0.946 | 0.955 | 0.950 |
| | 12 | 0.954 | 0.952 | 0.948 |
| | 16 | 0.952 | 0.954 | 0.950 |
| | 20 | 0.945 | 0.948 | 0.950 |
| 0.50 | 4 | 0.950 | 0.952 | 0.949 |
| | 8 | 0.944 | 0.952 | 0.948 |
| | 12 | 0.954 | 0.953 | 0.944 |
| | 16 | 0.946 | 0.950 | 0.949 |
| | 20 | 0.952 | 0.948 | 0.949 |
| 0.75 | 4 | 0.955 | 0.952 | 0.953 |
| | 8 | 0.953 | 0.948 | 0.949 |
| | 12 | 0.950 | 0.946 | 0.947 |
| | 16 | 0.948 | 0.946 | 0.951 |
| | 20 | 0.956 | 0.946 | 0.951 |

Finally observing that the posterior of $\eta$ in (6.3) is log-concave, the $100(1-\alpha)\%$ HPD interval $[\eta_1, \eta_2]$ for $\eta$ is obtained by solving

$$\eta_1^{n-2}(\eta_1^2 + \frac{S_{22.1}}{S_{11}})^{-(n-3/2)} = \eta_2^{n-2}(\eta_2^2 + \frac{S_{22.1}}{S_{11}})^{-(n-3/2)},$$

where

$$c \int_{\eta_1}^{\eta_2} \eta^{n-2}(\eta^2 + \frac{S_{22.1}}{S_{11}})^{-(n-3/2)} \, d\eta = 1 - \alpha,$$

$c$ being the normalizing constant.

Now we evaluate the frequentist coverage probability by investigating the HPD credible interval of the marginal posterior densities of $\beta$, $\theta$ and $\eta$ under our probability matching prior for several $\rho$ and $n$. That is to say, the frequentist coverage of a $100(1-\alpha)\%$ HPD interval should be close to $1-\alpha$. This is done numerically. The results were fairly insensitive to the choice of $\sigma_1$ and $\sigma_2$. Table 1 gives numerical values of the frequentist coverage probabilites of 95% HPD intervals for $\beta$, $\theta$ and $\eta$ for $\sigma_1 = \sigma_2 = 1$.

The computation of these numerical values is based on simulation. In particular, for fixed $(\mu_1, \mu_2, \sigma_1^2, \sigma_2^2, \rho)$ and $n$, we take 5,000 independent random samples of $(\boldsymbol{X}_1, \boldsymbol{X}_2)$ from the bivariate normal model. In our simulation study, we take $\mu_1 = \mu_2 = 0$ without loss of generality. Under the prior $\pi$, the frequentist coverage probability can be estimated by the relative frequency of HPD intervals containing the true parameter value. An inspection of Table 1 reveals that the agreement between the frequentist and posterior coverage probabilities of HPD intervals is quite good for the probability matching prior even if $n$ is small.

## 7. Summary

The paper considers several probability matching criteria, and develops a prior that meets all the matching criteria individually for several parameters of the bivariate normal distribution including the regression coefficient and the generalized variance. Future work will address development of matching priors when the parameter of interest is the correlation coefficient. Possible multivariate extensions will also be considered.

**Acknowledgments.** The paper has benefitted much from the comments of Bertrand Clarke and a reviewer.